\magnification=\magstep1

\newfam\msbfam 
\font\tenmsb=msbm10 \textfont\msbfam=\tenmsb  
\font\sevenmsb=msbm7 \scriptfont\msbfam=\sevenmsb

\def\hexnumber#1{\ifcase#1 0\or1\or2\or3\or4\or5\or6\or7\or8\or9\or
 A\or B\or C\or D\or E\or F\fi}
\edef\msbhx{\hexnumber\msbfam}
\mathchardef\ltimes="2\msbhx6E

\newfam\euffam
\font\teneuf=eufm10 \textfont\euffam=\teneuf 
\font\seveneuf=eufm7 \scriptfont\euffam=\seveneuf 
\def\frak#1{{\fam\euffam\relax#1}}
\def\Lc{{\frak c}}
\def\Lg{{\frak g}}
\def\Lk{{\frak k}}
\def\Ll{{\frak l}}
\def\Lm{{\frak m}}
\def\Ln{{\frak n}}
\def\Lr{{\frak r}}
\def\Ls{{\frak s}}
\def\Lz{{\frak z}}

\def\Lgl{{\Lg\Ll}}

\def\GL{{\rm GL}}
\def\SL{{\rm SL}}

\def\Re{{\rm Re}}
\def\End{{\rm End}}
\def\Ad{{\rm Ad}}
\def\ad{{\rm ad}}
\def\Stab{{\rm Stab}}

\def\barG{{\overline{G}}}
\def\barn{{\overline{n}}}
\def\barg{{\overline{g}}}
\def\barpi{{\overline{\pi}}}
\def\barLg{{\overline{\Lg}}}
\def\bargi{{\overline{g_i}}}

\def\tildeS{{\widetilde{S}}}
\def\tildeLs{{\widetilde{\Ls}}}
\def\tildeX{{\widetilde{X}}}
\def\tildeY{{\widetilde{Y}}}
\def\tildeT{{\widetilde{T}}}

\def\scrh{{\cal H}}

\newfam\msbfam 
\font\tenmsb=msbm10 \textfont\msbfam=\tenmsb  
\font\sevenmsb=msbm7 \scriptfont\msbfam=\sevenmsb

\def\Bbb#1{{\fam\msbfam\relax#1}}
\def\C{{\Bbb C}}
\def\R{{\Bbb R}}

\centerline{\bf Decay to zero of matrix coefficients at Adjoint infinity}

\medskip
\centerline{by}
\smallskip
\centerline{Scot Adams}

\medskip\bigskip\medskip\bigskip\noindent
{\bf I. Introduction}

\medskip\bigskip
The main theorems below are Theorem 9 and Theorem 11. As far
as I know, Theorem~9 represents a slight improvement over what
currently appears in the literature, and gives a fairly easy proof of
Theorem 11 which is due to R.~Howe and C.~Moore [4]. 
In the semisimple case, the Howe-Moore result follows from [6] or [7].

The proofs appearing here are relatively elementary and some readers
will recognize the influence of R.~Ellis and M.~Nerurkar [2].
What may be less evident is the connection to
N.~Kowalsky's work [5] on dynamics on Lorentz manifolds:
Lemma 1 below is a unitary analogue of the fact that
an expansive Adjoint action can result in
much of the Lie algebra being lightlike.
In Hilbert space, the situation is even nicer:
an isotropic vector must equal zero,
so we get that the Lie algebra of the stabilizer contains
all ``Kowalsky'' vectors for the Adjoint sequence,
and not just a codimension one subspace.
Moreover, much of the rest of proof of Theorem 9
also uses ideas that were originally developed to describe effectively
the collection of simply connected Lie groups admitting
an orbit nonproper action on a connected Lorentz manifold.
Thus the debt to Kowalsky is significant.

In this note, we assume some familiarity with Lie theory
and basic unitary representation theory.
The exposition should otherwise be self-contained.

Please send comments, suggestions or questions to adams@math.umn.edu.

Throughout, ``Lie group'' means ``$C^\infty$ real Lie group'',
``Lie algebra'' means ``real Lie algebra''
 and ``Hilbert space'' means ``complex Hilbert space''.
Lie groups are denoted by capital Roman letters and, for any Lie group,
its Lie algebra is denoted, without comment, by the corresponding
small letter in the fraktur font.
The sesquilinear form on a Hilbert space is denoted
$\langle\,\,\cdot\,\,,\,\,\cdot\,\,\rangle$.
Convergence in the weak topology on a Hilbert
space is denoted with $\rightharpoonup$, so
$v_i\rightharpoonup v$ in a Hilbert space $\scrh$ means:
for all $w\in\scrh$, $\langle v_i,w\rangle\to\langle v,w\rangle$.
The group of bounded operators on $\scrh$ is denoted $B(\scrh)$.
The unitary group of a Hilbert space~$\scrh$ is denoted $U(\scrh)$.
The weak-operator and strong-operator topologies agree on $U(\scrh)$,
and $U(\scrh)$ is given this topology.
Unitary representations are always assumed
continuous in the weak-operator topology
(or, equivalently, in the strong-operator topology).

I do not plan to publish this note.

\medskip\bigskip\medskip\bigskip\noindent
{\bf II. Preliminary dynamical results}

\medskip\bigskip
Let $G$ be a Lie group acting by
a unitary representation on a Hilbert space~$\scrh$.
Let $w\in\scrh$ and let $\Lm$ be the Lie algebra of $\Stab_G(w)$.

\medskip\bigskip\noindent
{\bf Lemma 1:} {\sl Let $U_i\to0$ in $\Lg$.
Let $g_i$ be a sequence in $G$.
Assume $(\Ad\,g_i)U_i\to X$ in $\Lg$.
Let $v\in\scrh$ and assume that $g_iv\rightharpoonup w$ in $\scrh$.
Then $X\in\Lm$.}

\medskip\noindent
{\it Proof:} Fix $t\in\R$ and let $h:=\exp(tX)$. We wish to show that $hw=w$.

As $\exp(tU_i)\to1_G$ in $G$,
we get $(\exp(tU_i))v\to v$ in $\scrh$.
So, for all $x\in\scrh$,
since $\{g_i^{-1}x\}$ is norm-bounded in $\scrh$,
Cauchy-Schwartz yields:
$\langle[(\exp(tU_i))v]-v,g_i^{-1}x\rangle\to0$; then
$$\langle g_i(\exp(tU_i))v,x\rangle-\langle g_iv,x\rangle\to0.$$
Moreover, as $g_iv\rightharpoonup w$, we see, for all $x\in\scrh$, that
$$\langle g_iv,x\rangle-\langle w,x\rangle\to0.$$
Adding the last two displayed limits, for all $x\in\scrh$,
we have $\langle g_i(\exp(tU_i))v,x\rangle-\langle w,x\rangle\to0$.
That is, $g_i(\exp(tU_i))v\rightharpoonup w$ in $\scrh$.

For all $i$, let $X_i:=(\Ad\,g_i)U_i$ and let $h_i:=\exp(tX_i)$.
Then, for all $i$, we have $g_i(\exp(tU_i))g_i^{-1}=h_i$,
so $h_ig_i=g_i(\exp(tU_i))$.
Then $h_ig_iv=g_i(\exp(tU_i))v\rightharpoonup w$ in $\scrh$.

We have $X_i\to X$ in $\Lg$, and so $h_i\to h$ in $G$.
For all $x\in\scrh$, since $h_i^{-1}x\to h^{-1}x$
and since $\{g_iv\}$ is norm-bounded in $\scrh$,
Cauchy-Schwartz gives:
$\langle g_iv,h_i^{-1}x-h^{-1}x\rangle\to0$;
then
$$\langle h_ig_iv,x\rangle-\langle hg_iv,x\rangle\to0.$$
Since $g_iv\rightharpoonup w$, it follows, for all $x\in\scrh$,
that $\langle g_iv,h^{-1}x\rangle-\langle w,h^{-1}x\rangle\to0$;
then
$$\langle hg_iv,x\rangle-\langle hw,x\rangle\to0.$$
Adding the last two displayed limits, for all $x\in\scrh$, we have
$\langle h_ig_iv,x\rangle-\langle hw,x\rangle\to0$.
That is, $h_ig_iv\rightharpoonup hw$.
So, recalling that $h_ig_iv\rightharpoonup w$,
we get $hw=w$. QED

\medskip\bigskip\noindent
{\bf Corollary 2:} {\sl For all $T\in\Lm$,
we have $(\Lc_\Lg(T))\cap((\ad\,T)\Lg)\subseteq\Lm$.}

\medskip\noindent
{\it Proof:}
Let $X\in(\Lc_\Lg(T))\cap((\ad\,T)\Lg)$.
We wish to show that $X\in\Lm$.

We have $(\ad\,T)X=0$ and $X\in(\ad\,T)\Lg$.
Choose $S\in\Lg$ such that $(\ad\,T)S=X$.
Let $r_i$ be a sequence of nonzero real numbers such that $r_i\to+\infty$.
For all $i$, let $g_i:=\exp(r_iT)$ and let $U_i:=S/r_i$.
Then $U_i\to0$ in $\Lg$.
Since $T\in\Lm$, it follows, for all $i$, that $g_iw=w$.

Since $(\ad\,T)S=X$ and $(\ad\,T)X=0$,
it follows, for all $i$, that $(\Ad\,g_i)S=S+r_iX$.
Then $(\Ad\,g_i)U_i=(S/r_i)+X\to X$ in $\Lg$,
so, by Lemma 1 (with $v:=w$), we are done. QED

\medskip\bigskip\medskip\bigskip\noindent
{\bf III. More preliminaries}

\medskip\bigskip\noindent
{\bf Remark 3:} {\sl Let $\Lg$ be a Lie algebra.
Let $U_i\to0$ in $\Lg$ and let $X\in\Lg$.
Let $g_i$ be a sequence in~$G$.
Assume that $(\Ad\,g_i)U_i\to X$ in $\Lg$.
Then $\ad\,X:\Lg\to\Lg$ is nilpotent.}

\medskip\noindent
{\it Proof:}
Let $F_i\subseteq\C$ be the set of eigenvalues of
$\ad\,U_i:\Lg\to\Lg$.
Then, as $U_i\to0$,
it follows that $F_i\to\{0\}$
in the topological space of finite subsets of $\C$.
For all $i$, in $\End(\Lg)$, we have
$$\ad((\Ad\,g_i)U_i)\qquad=\qquad(\Ad\,g_i)\,(\ad\,U_i)\,(\Ad\,g_i)^{-1},$$
so $F_i$ is also the set of eigenvalues of $\ad((\Ad\,g_i)U_i):\Lg\to\Lg$.
Passing to the limit,
$\{0\}$ is the set of eigenvalues of $\ad\,X:\Lg\to\Lg$. QED

\medskip\bigskip\noindent
{\bf Remark 4:} {\sl Let $V$ be a finite-dimensional real vector space.

\itemitem{(i)} Let $T_i$ be a sequence in $\GL(V)$.
Assume that $\{T_i\}$ is not precompact in $\GL(V)$.
Then either $\{T_i\}$ or $\{T_i^{-1}\}$
is not precompact in $\End(V)$.
\itemitem{(ii)} Let $E_i$ be a sequence in $\End(V)$.
Assume that $\{E_i\}$ is not precompact in $\End(V)$.
Then there exists $v\in V$ such that $E_i(v)$ is not precompact
in $V$.

}

\medskip\noindent
{\it Proof:} {\it Proof of (i):}
For all $i$, let
$U_i:= T_i\oplus T_i^{-1}\in\GL(V\oplus V)$,
so $U_i$ is defined by $U_i(v,w)=(T_iv,T_i^{-1}w)$.
Then $\{U_i\}$ is not precompact in~$\GL(V\oplus V)$.
For all $i$, we have $\det(U_i)=1$,
so $\{U_i\}$ is not precompact in~$\SL(V\oplus V)$.
So, since $\SL(V\oplus V)$ is closed in $\End(V\oplus V)$,
we see that $\{U_i\}$ is not precompact in $\End(V\oplus V)$.
Then, as $\{U_i\}=\{T_i\oplus T_i^{-1}\}$,
the result follows.
{\it End of proof of~(i).}

{\it Proof of (ii):}
Let $d:=\dim(V)$.
Let $v_1,\ldots,v_d$ be a basis for $V$.
Let $\psi_1,\ldots,\psi_d:V\to\R$ be a basis for the dual space of $V$.
For all $i$, let $M_i:=[\psi_k(E_i(v_j))]_{jk}$ be
the matrix of $E_i$ with respect to the two bases.
Then $\{M_i\}_i$ is not precompact in $\R^{d\times d}$.
Choose $j,k\in\{1,\ldots,d\}$ such that
$\{\psi_k(E_i(v_j))\}_i$ is not precompact in $\R$.
Then $\{E_i(v_j)\}_i$ is not precompact in $V$.
Let $v:=v_j$.
{\it End of proof of (ii).} QED

\medskip\bigskip
Let $\Lg$ be a Lie algebra and let $X,Y,T\in\Lg$.
We will say that
$(X,Y,T)$ is a {\bf standard triple} for $\Lg$ if
$[T,X]=2X$, $[T,Y]=-2Y$ and $[X,Y]=T\ne0$.

Let $G$ be a Lie group and let $X,Y,T\in\Lg$.
For all $u\in\R$, let
$$n(u):=\exp(uX),\qquad \barn(u):=\exp(uY),\qquad a(u):=\exp(uT).$$
We will say that $(X,Y,T)$ is a {\bf basic triple} for $G$ if,
for all $\tau\in\R$, for all $\delta\in\R\backslash\{0\}$,
$$\left[n\left({e^\tau-1\over\delta}\right)\right]
\quad\cdot\quad[\barn(\delta)]\quad\cdot\quad
\left[n\left({e^{-\tau}-1\over\delta}\right)\right]
\quad=\quad\left[\barn({e^{-\tau}\delta})\right]\quad\cdot\quad[a(\tau)].$$

\medskip\bigskip\noindent
{\bf Remark 5:} {\sl 
Let $G$ be a Lie group.
Then any standard triple for $\Lg$ is a basic triple for $G$.}

\medskip\noindent
{\it Proof:}
Let $(X,Y,T)$ be a standard triple for $\Lg$.
We wish to show that $(X,Y,T)$ is a basic triple for $G$.
Let $S:=\SL_2(\R)$ and let
$$X_0=\left[\matrix{
0  & 1 \cr
0  & 0 \cr
}\right],\qquad
Y_0=\left[\matrix{
0  & 0 \cr
1  & 0 \cr
}\right],\qquad
T_0=\left[\matrix{
1  & 0 \cr
0  &-1 \cr
}\right].$$
Then, for all $u\in\R$, we have
$$\exp(uX_0)=\left[\matrix{
1  & u \cr
0  & 1 \cr
}\right],\qquad
\exp(uY_0)=\left[\matrix{
1  & 0 \cr
u  & 1 \cr
}\right],\qquad
\exp(uT_0)=\left[\matrix{
e^u  &   0   \cr
0    & e^{-u}\cr
}\right],$$
so straightforward matrix computations prove that
$(X_0,Y_0,T_0)$ is a standard triple for~$\Ls$
and a basic triple for $S$.
Let $\tildeS$ be the universal cover of $S$
and let $p:\tildeS\to S$ be a covering homomorphism.
Then $dp:\tildeLs\to\Ls$ is an isomorphism of Lie algebras.
Let $\tildeX_0:=(dp)^{-1}(X_0)$,
let $\tildeY_0:=(dp)^{-1}(Y_0)$ and
let $\tildeT_0:=(dp)^{-1}(T_0)$.
Then $(\tildeX_0,\tildeY_0,\tildeT_0)$ is a standard
triple for $\tildeLs$ and, by uniquness of liftings across covering maps,
$(\tildeX_0,\tildeY_0,\tildeT_0)$ is a basic triple for $\tildeS$.

Since $(\tildeX_0,\tildeY_0,\tildeT_0)$ is a standard triple for $\tildeLs$
and $(X,Y,T)$ is a standard triple for $\Lg$,
there is a homomorphism of Lie algebras
$\phi:\tildeLs\to\Lg$ such that
$\phi(\tilde X_0)=X$, $\phi(\tilde Y_0)=Y$ and $\phi(\tilde T_0)=T$.
By the Monodromy Theorem (see  Theorem 2.7.5, on p.~71 of [8]),
there is a Lie group homomorphism $f:\tildeS\to G$ such that $df=\phi$.
Then, by naturality of the exponential map, we see that $(X,Y,T)$ is
a basic triple for $G$.
QED

\medskip\bigskip\noindent
{\bf Proposition 6:} {\sl Let $S$ be a connected simple Lie group
acting on a Hilbert space $\scrh$ by a unitary representation.
Let $w\in\scrh$.
Let $\Lm$ be the Lie algebra of $\Stab_S(w)$.
Assume $\exists X\in\Lm\backslash\{0\}$ such that $\ad\,X:\Ln\to\Ln$ is nilpotent.
Then $S$ fixes $w$.}

\medskip\noindent
{\it Proof:}
By Jacobson-Morozov
(see, {\it e.g.}, Theorem 7.4 of Chapter IX on p.~432 of [3]),
choose $Y,T\in\Ls$ such that 
$(X,Y,T)$ is a standard triple for $\Ls$.
Then $T\ne0$.
For all $u\in\R$, let
$$n(u):=\exp(uX),\qquad
\barn(u):=\exp(uY),\qquad
a(u):=\exp(uT).$$
Let $A:=\{a(u)\}_{u\in\R}$ be the image of $a:\R\to S$.

{\it Claim $\alpha$:} $w$ is $A$-invariant.
{\it Proof of Claim $\alpha$:}
Let $N:=\{n(u)\}_{u\in\R}$ be the image of~$n:\R\to S$.
Since $X\in\Lm$, we see that $w$ is $N$-invariant.

Define $f:S\to\R$ by $f(s)=\Re(\langle sw,w\rangle)$.
For all $s\in S$, and all $p,q\in N$,
we have $\langle psq w,w\rangle=
\langle s(qw),p^{-1}w\rangle=
\langle sw,w\rangle$,
so $f(psq)=f(g)$.
That is, $f$ is bi-invariant under~$N$.

Let $a_0\in A$. We wish to show that $a_0w=w$.
We have $f(1_S)=\|w\|^2$ and
$$2(f(a_0))=-\|a_0w-w\|^2+\|a_0w\|^2+\|w\|^2=-\|a_0w-w\|^2+2\|w\|^2,$$
so $2(f(a_0))-2(f(1_S))=-\|a_0w-w\|^2$.
It therefore suffices to prove $f(a_0)=f(1_S)$.

Fix $\tau\in\R$ such that $a_0=a(\tau)$.
Let $\delta_i$ be a sequence of nonzero real numbers
such that $\delta_i\to0$ in $\R$.
For all $i$, let
$$n_i:=n\left({e^\tau-1\over\delta_i}\right),\qquad
\barn_i:=\barn(\delta_i),\qquad
n'_i:=n\left({e^{-\tau}-1\over\delta_i}\right).$$
By Remark 5, since $(X,Y,T)$ is a standard triple for $\Ls$,
$(X,Y,T)$ is a basic triple for $S$.
Then, for all $i$,
we have $n_i\barn_in'_i=[\barn(e^{-\tau}\delta_i)]\cdot a_0$,
so, since $\delta_i\to0$ in $\R$,
we get $n_i\barn_in'_i\to a_0$ in $S$.
Moreover, $\barn_i=\barn(\delta_i)\to\barn(0)=1_S$.
By bi-invariance of $f$ under $N$, for all $i$,
we have $f(n_i\barn_i n'_i)=f(\barn_i)$,
so, passing to the limit, $f(a_0)=f(1_S)$.
{\it End of proof of Claim $\alpha$.}

For all $\lambda\in\R$, let $\Ls_\lambda:=\{A\in\Ls\,|\,(\ad\,T)A=\lambda A\}$.

{\it Claim $\beta$:} For all $\lambda\in\R\backslash\{0\}$,
we have $\Ls_\lambda\subseteq\Lm$.
{\it Proof of Claim $\beta$:}
Fix $\lambda\in\R\backslash\{0\}$ and $A\in\Ls_\lambda$.
We wish to show that $A\in\Lm$.

Fix a sequence $r_i$ in $\R$ such that $r_i\lambda\to+\infty$.
For all $i$, let $g_i:=a(r_i)=\exp(r_iT)$;
then, as $A\in\Ls_\lambda$,
we get $(\ad\,T)A=\lambda A$, so $(\Ad\,g_i)A=e^{r_i\lambda}A$,
so $(\Ad\,g_i)(e^{-r_i\lambda}A)=A$.
Moreover, by Claim $\alpha$, we have $g_iw=w$.
Then, by Lemma 1
(with $v:=w$, with $U_i$~replaced by~$e^{-r_i\lambda}A$
and with $X$ replaced by $A$), we are done.
{\it End of proof of Claim~$\beta$.}

Let $\Ls_*:=(\ad\,T)\Ls$ and
let $\Ls'$ be the Lie subalgebra of $\Ls$ generated by the set $\Ls_*$.
By, {\it e.g.}, Corollary 4.11.7 on p.~130 of [1],
$\ad\,T:\Ls\to\Ls$ is real diagonalizable.
Then $\displaystyle{\Ls_*=
\sum_{\lambda\in\R\backslash\{0\}}\,\Ls_\lambda}$ and
$\Ls=\Ls_0+\Ls_*$ and,
from Claim $\beta$, we get $\Ls_*\subseteq\Lm$.
Then $\Ls'\subseteq\Lm$.

The Jacobi identity implies that $\Lc_\Ls(T)$ normalizes the set $(\ad\,T)\Ls$,
{\it i.e.}, that $\Ls_0$ normalizes~$\Ls_*$.
Then $\Ls_0$ normalizes $\Ls'$.
Moreover, $\Ls_*\subseteq\Ls'$, so $\Ls_*$ normalizes~$\Ls'$.
Then $\Ls_0+\Ls_*$ normalizes~$\Ls'$, {\it i.e.}, $\Ls$ normalizes~$\Ls'$.
That is, $\Ls'$ is an ideal in $\Ls$.
Since $T\ne0$, we have $\ad_\Ls(T)\ne0$, so $\Ls_*\ne\{0\}$, so $\Ls'\ne\{0\}$.
Then, by simplicity of $\Ls$, it follows that $\Ls'=\Ls$.
Then $\Ls=\Ls'\subseteq\Lm$.
That is, $S$ fixes $w$.
QED

\medskip\bigskip
The proof of Claim $\alpha$ above is an algebraic
argument derived from the geometry
of the proof of Theorem 2.4.2 on p.~29 of [9].

The following is Lemma 4.15.7 on p.~144 of [1]:

\medskip\noindent
{\bf Lemma 7:} {\sl Let $\Lg$ be a Lie algebra with no simple direct summand
and let $\Ln$ be the nilradical of $\Lg$.
Then $\Lc_\Lg(\Ln)=\Lz(\Ln)$.}

\medskip\noindent
{\it Proof:}
Let $\Ls$ be a Levi factor of $\Lg$ and let $\Lr$ be the solvable radical of $\Lg$.
Then, because $\Lc_\Lg(\Ln)$ is an ideal in $\Lg$, it follows
(by, {\it e.g.}, Lemma 4.10.19 on p.~119 of [1]) that
$$\Lc_\Lg(\Ln)\qquad=\qquad[(\Lc_\Lg(\Ln))\cap\Ls]\qquad+\qquad[(\Lc_\Lg(\Ln))\cap\Lr],$$
{\it i.e.}, that $\Lc_\Lg(\Ln)=(\Lc_\Ls(\Ln))+(\Lc_\Lr(\Ln))$.

By (ii) of Theorem 3.8.3 on p.~206 of [8],
we have $[\Lg,\Lr]\subseteq\Ln$.
It follows that $[\Ls,\Lr]\subseteq\Ln$,
{\it i.e.}, that $\ad:\Ls\to\Lgl(\Lr/\Ln)$ is zero.
On the other hand, as $\Lg$ has no simple direct summand,
it follows that $\ad:\Ls\to\Lgl(\Lr)$ is faithful.
Therefore, by complete reducibility of $\ad:\Ls\to\Lgl(\Lr)$,
$\ad:\Ls\to\Lgl(\Ln)$ is faithful.
That is, $\Lc_\Ls(\Ln)=\{0\}$.
Then $\Lc_\Lg(\Ln)=(\Lc_\Ls(\Ln))+(\Lc_\Lr(\Ln))=\Lc_\Lr(\Ln)$.

We have $[\Lc_\Lr(\Ln),\Lc_\Lr(\Ln)]\subseteq[\Lg,\Lr]\subseteq\Ln$
and $[\Lc_\Lr(\Ln),\Ln]=\{0\}$.
Then $\Lc_\Lr(\Ln)$ is a nilpotent ideal of $\Lg$,
so $\Lc_\Lr(\Ln)\subseteq\Ln$.
Then $\Lc_\Lg(\Ln)=\Lc_\Lr(\Ln)\subseteq\Ln$,
so $\Lc_\Lg(\Ln)=(\Lc_\Lg(\Ln))\cap\Ln=\Lz(\Ln)$.
QED

\medskip\bigskip\noindent
{\bf Lemma 8:} {\sl Let $\Lg$ be a Lie algebra with no simple direct summand.
Let $\Ln$ be the nilradical of $\Lg$.
Let $\Lm$ be a vector subspace of $\Lg$.
Assume $\forall T\in\Lm$, $(\Lc_\Lg(T))\cap((\ad\,T)\Lg)\subseteq\Lm$.
Assume $\exists X\in\Lm\backslash\{0\}$ such that $\ad\,X:\Ln\to\Ln$ is nilpotent.
Then $\Lm\cap(\Lz(\Ln))$ is a nonzero ideal of $\Lg$.}

\medskip\noindent
{\it Proof:}
As $\Lz(\Ln)$ is an Abelian ideal of $\Lg$, for all $T\in\Lz(\Ln)$,
we get $(\ad\,T)\Lg\subseteq(\Lc_\Lg(T))\cap(\Lz(\Ln))$.
Then, for all $T\in\Lm\cap(\Lz(\Ln))$, we have
$(\ad\,T)\Lg\subseteq(\Lc_\Lg(T))\cap((\ad\,T)\Lg)\cap(\Lz(\Ln))
\subseteq\Lm\cap(\Lz(\Ln))$.
That is, $\Lm\cap(\Lz(\Ln))$ is an ideal of $\Lg$.
It remains to show that $\Lm\cap(\Lz(\Ln))\ne\{0\}$.

Define the descending central series of $\Ln$ by
$\Ln^{(1)}:=\Ln$ and $\Ln^{(i+1)}:=[\Ln,\Ln^{(i)}]$.
Fix an integer $l\ge1$ such that $\Ln^{(l)}=\{0\}$.
Say a sequence $P_0,P_1,P_2,\ldots$ in $\Lg$ is {\bf good} if
\itemitem{(1)} $P_0\in\Lm\backslash\{0\}$;
\itemitem{(2)} $\ad\,P_0:\Ln\to\Ln$ is nilpotent; and
\itemitem{(3)} for all integers $i\ge0$, we have
$P_{i+1}\in(\Lc_\Ln(P_i))\cap((\ad\,P_i)\Ln)$.

\noindent
By assumption, $X,0,0,0,\ldots$ is good.
For any good sequence $P_i$, by induction, for all integers $i\ge1$,
we have $P_i\in\Ln^{(i)}$; in particular, $P_l=0$.

Fix a good sequence $P_i$ such that
$k:=\max\{i\,|\,P_i\ne0\}$ is as large as possible.
Then $(\Lc_\Ln(P_k))\cap((\ad\,P_k)\Ln)=\{0\}$.
If $k=0$, then, by (2), $\ad_\Ln(P_k)$ is nilpotent.
If $k\ge1$, then $P_k\in\Ln^{(k)}\subseteq\Ln$, and so, again,
$\ad_\Ln(P_k)$ is nilpotent.

Then $(\Lc_\Ln(P_k))\cap((\ad\,P_k)\Ln)$
is intersection of the kernel and image of the nilpotent map
$\ad_\Ln(P_k)$.
However, $(\Lc_\Ln(P_k))\cap((\ad\,P_k)\Ln)=\{0\}$,
while the intersection of the kernel and image of a nonzero
nilpotent endomorphism is never $\{0\}$.
Then $\ad_\Ln(P_k)=0$, {\it i.e.}, $P_k\in\Lc_\Lg(\Ln)$.
So, as $\Lg$ has no simple direct summand, by Lemma 7,
we get $P_k\in\Lz(\Ln)$.

For all integers $i\ge0$, we have
$P_{i+1}\in(\Lc_\Ln(P_i))\cap((\ad\,P_i)\Ln)
\subseteq(\Lc_\Lg(P_i))\cap((\ad\,P_i)\Lg)$.
So, by induction, for all integers $i\ge0$, we have $P_i\in\Lm$.
Then $0\ne P_k\in\Lm\cap(\Lz(\Ln))$.
QED

\medskip\bigskip\medskip\bigskip\noindent
{\bf IV. Decay to zero at Adjoint infinity for connected Lie groups}

\medskip\bigskip\noindent
{\bf Theorem 9:} {\sl Let $G$ be a connected Lie group
acting on a Hilbert space $\scrh$ by
a unitary representation $\pi:G\to U(\scrh)$.
Assume that

\noindent
\hskip.1in$(*)$ no nonzero vector of $\scrh$ is fixed by a nontrivial
normal connected subgroup of $G$.

\noindent
Let $g_i$ be a sequence in $G$ and
assume that $\Ad_\Lg(g_i)$ leaves compact subsets of $\GL(\Lg)$.
Then $\pi(g_i)\to0$ in the weak-operator topology on $B(\scrh)$.}

\medskip\bigskip\noindent
{\bf Notes:}
\item{1.} Condition $(*)$ is satisfied if $G$ is simple and $\scrh$
admits no $G$-invariant vectors.
\item{2.} Condition $(*)$ is satisfied if $\pi$ is faithful and irreducible.
(The set of vectors fixed by a
normal subgroup is a $G$-invariant subspace.
By irreducibility, if such a subspace were nonzero,
it would equal $\scrh$.
Then, by faithfulness, the subgroup would be trivial.)
\item{3.} I am not sure whether, in Theorem 9,
we need that $G$ is connected.
\item{4.} One may summarize Theorem 9 as asserting:
If a unitary representation of a connected Lie group satisfies $(*)$,
then its matrix coefficients decay to zero at ``$\Ad$-infinity''.

\medskip\bigskip\noindent
{\it Proof of Theorem 9:}
Let $B(\scrh)$ have the weak-operator topology.
Assume, for a contradiction, that $\pi(g_i)\not\to0$ in $B(\scrh)$.

Passing to a subsequence, assume that $\{\pi(g_i)\}$ is
bounded away from $0$ in $B(\scrh)$.
For all $i$, $\pi(g_i^{-1})=(\pi(g_i))^*$,
so, because $L\mapsto L^*:B(\scrh)\to B(\scrh)$ is continuous,
we see that $\{\pi(g_i^{-1})\}$ is bounded away from $0$ in $B(\scrh)$, as well.
By (i) of Remark 4, 
either $\{\Ad_\Lg(g_i)\}$ or $\{\Ad_\Lg(g_i^{-1})\}$
is not precompact in $\End(\Lg)$.
By, if necessary, replacing each $g_i$ with $g_i^{-1}$,
assume that  $\{\Ad_{\Lg}(g_i)\}$ is not precompact in $\End(\Lg)$.
Passing to a subsequence, assume that $\Ad_\Lg(g_i)$ leaves
compact subsets of $\End(\Lg)$.
Fix $v\in\scrh$ such that $g_iv\not\rightharpoonup0$ in $\scrh$.
Passing to a subsequence, assume that $g_iv\rightharpoonup w\ne0$ in $\scrh$.

Choose normal connected Lie subgroups $G'$ and $G''$ of $G$
such that $\Lg=\Lg'+\Lg''$,
such that $\Lg'\cap\Lg''=\{0\}$,
such that $\Lg'$ is semisimple or trivial and
such that $\Lg''$ has no simple direct summand.
It follows either that
$\{\Ad_{\Lg'}(g_i)\}$ is not precompact in~$\End(\Lg')$
or that
$\{\Ad_{\Lg''}(g_i)\}$ is not precompact in $\End(\Lg'')$.

{\it Case A: $\{\Ad_{\Lg'}(g_i)\}$ is not precompact in~$\End(\Lg')$.}
Choose a simple normal connected Lie subgroup $S$ of $G'$
such that $\{\Ad_\Ls(g_i)\}$ is not precompact in $\End(\Ls)$.
By (ii) of Remark~4,
choose $U\in\Ls$ such that $\{(\Ad\,g_i)U\}$ is not precompact in $\Ls$.

Passing to a subsequence, choose $t_i\to0$ in $\R$
such that $t_i(\Ad\,g_i)U\to X\ne0$ in $\Ls$.
By Remark 3 (with $\Lg$ replaced by $\Ls$ and $U_i$ replaced by $t_iU$),
$\ad\,X:\Ls\to\Ls$ is nilpotent.
Let $\Lm$ be the Lie algebra of $\Stab_S(w)$.
By Lemma 1 (with $G$ replaced by $S$ and $U_i$ replaced by~$t_iU$),
we see that $X\in\Lm$.
Then, by Proposition 6, $S$ fixes $w$.
However, $S$~is normal in~$G'$, and therefore in $G$,
contradicting $(*)$.
{\it End of Case A.}

{\it Case B: $\{\Ad_{\Lg''}(g_i)\}$ is not precompact in $\End(\Lg'')$.}
By (ii) of Remark~4, choose $U\in\Lg''$ such that
$\{(\Ad\,g_i)U\}$ is not precompact in $\Lg''$.

Passing to a subsequence, choose $t_i\to0$ in $\R$
such that $t_i(\Ad\,g_i)U\to X\ne0$ in~$\Lg''$.
By Remark 3 (with $\Lg$ replaced by $\Lg''$ and $U_i$ replaced by $t_iU$),
$\ad\,X:\Lg''\to\Lg''$ is nilpotent.
Let $\Lm$ be the Lie algebra of $\Stab_{G''}(w)$.
By Lemma 1 (with $G$~replaced by $G''$ and $U_i$~replaced by $t_iU$),
we get $X\in\Lm$.
By Corollary~2 (with $G$~replaced by~$G''$),
for all $T\in\Lm$, we have
$(\Lc_{\Lg''}(T))\cap((\ad\,T)\Lg'')\subseteq\Lm$.
Let $\Ln$ be the nilradical of~$\Lg''$.
Then $\ad\,X:\Ln\to\Ln$ is nilpotent.
Then, by Lemma 8 (with $\Lg$ replaced by $\Lg''$),
$\Lm\cap(\Lz(\Ln))$ is a nonzero ideal of~$\Lg''$,
and therefore of $\Lg$.
Then the connected Lie subgroup of $G$ corresponding to $\Lm\cap(\Lz(\Ln))$
is a nontrivial normal connected subgroup of $G$ fixing $w$,
contradicting $(*)$.
{\it End of Case~B.} QED

\medskip\bigskip
Thanks to D.~Witte Morris for help in developing the following example:

\medskip\noindent
{\bf Example 10:} {\sl In the statement of Theorem 9,
$(*)$ cannot be replaced by

\noindent
\hskip.1in$(*')$ no nonzero vector of $\scrh$ is fixed by a noncompact
normal connected subgroup of $G$

\noindent
even under the assumption that $\pi$ is faithful.}

\medskip\noindent
{\it Proof:}
Let $H$ be the $3$-dimensional Heisenberg group,
and let $Z$ denote the center of $H$.
Let $\phi:H\to\R^2$ be a surjective homomorphism
whose kernel is $Z$,
and let $D$ be a nontrivial discrete subgroup of $Z$.
Let $G:=H/D$ and
let $p:H\to G$ be the canonical homomorphism.
Let $\psi:G\to\R^2$ be defined by $\psi(p(g))=\phi(g)$.

Let $\scrh'$ be a Hilbert space and let
$\rho:\R^2\to U(\scrh')$ be a unitary representation of~$\R^2$
with no nonzero invariant vectors
such that not all matrix coefficients decay to zero at infinity.
Let $\pi':=\rho\circ\psi:G\to U(\scrh')$.
Let $\scrh''$ be another Hilbert space and
fix a faithful unitary representation
$\pi'':G\to U(\scrh'')$.

Let $\scrh:=\scrh'\oplus\scrh''$.
Then $\pi:=\pi'\oplus\pi'':G\to U(\scrh)$
is faithful and satisfies $(*')$,
but does not enjoy the property that all
matrix coefficients decay to zero at $\Ad$-infinity.
QED

\medskip\bigskip\medskip\bigskip\noindent
{\bf V. Decay to zero at projective infinity for algebraic groups} 

\medskip\bigskip\noindent
{\bf Theorem 11:} {\sl Let $G$ be the connected real points of
a linear algebraic $\R$-group.
Let $\pi:G\to U(\scrh)$ be an irreducible unitary representation
on a Hilbert space $\scrh$.
Then $\pi(g)\to0$ in the weak-operator topology on $B(\scrh)$,
as $g$ leaves compact subsets of G modulo
the projective kernel of $\pi$.}

\medskip\noindent
{\it Proof:}
Let $K:=\ker(\pi)$.
Let $\barG:=G/K$.
Let $p:G\to\barG$ be the canonical homomorphism.
For all $g\in G$, let $\barg:=p(g)$.
Define $\barpi:\barG\to U(\scrh)$ by $\barpi(\barg)=\pi(g)$.

The Adjoint representation $\Ad:G\to\GL(\Lg)$ is algebraic
and $\Lk$ is an invariant subspace,
so $\Ad:G\to\GL(\Lg/\Lk)$ is algebraic and therefore
has closed image.
So, since $\barLg=\Lg/\Lk$ and since $\Ad_\barLg(\barG)=\Ad_\barLg(G)$,
we see that $\Ad_{\barLg}(\barG)$ is closed in $\GL(\barLg)$.
Then $\Ad:\barG\to\GL(\barLg)$ factors to a proper, injective
Lie group homomorphism $F:\barG/(Z(\barG))\to\GL(\barLg)$.

Assume that $g_i$ leaves compact subsets of G
modulo the projective kernel of $\pi$.
We wish to show that $\pi(g_i)\to0$
in the weak-operator topology on $U(\scrh)$.

The sequence $\bargi$ leaves compact subsets of $\barG$
modulo the projective kernel of $\barpi$.
By Schur's Lemma,
since $\barpi$ is faithful and irreducible,
its projective kernel is $Z(\barG)$.
So, since $F$ is proper,
$\Ad\,\bargi$ leaves compact subsets of $\GL(\barLg)$.
Then, by Theorem 9,
$\barpi(\bargi)\to0$ in the weak-operator topology on $U(\scrh)$.
For all $i$, $\barpi(\bargi)=\pi(g_i)$, so we are done.
QED

\medskip
The proof of Theorem 11 also works for any
connected Lie group that is locally isomorphic to a real linear
algebraic group.

\medskip\bigskip\medskip\bigskip
\centerline{\bf References}
\medskip\noindent
\item{[1]} S.~Adams: {\it Dynamics on Lorentz Manifolds}, World Scientific,
Singapore, 2001.
\item{[2]} R. Ellis and M. Nerurkar: {\it Weakly almost periodic flows},
Trans.~Amer.~Math.~Soc.~{\bf313} no.~1 (1989), 103--119.
\item{[3]} S.~Helgason: {\it Differential Geometry, Lie Groups and
Symmertic Spaces}, Academic Press, London, 1978.
\item{[4]} R.~Howe and C.~C.~Moore: {\it Asymptotic properties of unitary
representations.}, J.~Funct. Anal.~{\bf32} no.~1 (1979), 72--96.
\item{[5]} N.~Kowalsky: {\it oncompact simple automorphism groups of
Lorentz manifolds and other geometric manifolds},
Ann.~of Math.~(2) {\bf144} no.~3 (1996), 611--640.
\item{[6]} C.~C.~Moore: {\it Ergodicitiy of flows on homogeneous spaces},
Amer.~J.~Math.~{\bf88} (1966), 154--178.
\item{[7]} T.~Sherman: {\it A weight theory for unitary representations},
Canad.~J.~Math.~{\bf18} (1966), 159--168.
\item{[8]} V.~S.~Varadarajan: {\it Lie Groups, Lie Algebras, and Their
Representations}, Springer-Verlag, New York, 1984.
\item{[9]} R.~J.~Zimmer: {\it Ergodic Theory and Semisimple Groups}, Birkh\"auser,
Boston, 1984.

\end